\documentclass{elsart}
\usepackage{amsmath}
\usepackage{amscd}
\usepackage{amsfonts}

\journal{Topology and its Applications}

\newcommand{\sign}{\operatorname{sign}}

\newcommand{\R}{\mathbb{R}}

\newcommand{\N}{\mathbb{N}}
\newcommand{\G}{\Gamma}
\newcommand{\Hb}{\mathcal{H}}

\begin{document}
\begin{frontmatter}

\title{Group actions on Banach spaces and a geometric characterization of a-T-menability}
\author{Piotr W. Nowak}
\address{Department of Mathematics, Vanderbilt University, 1326 Stevenson Center, Nashville, TN 37240
USA.}

\begin{abstract}
We prove a geometric characterization of a-T-menability through
proper, affine, isometric actions on the Banach spaces $L_p[0,1]$
for $1<p<2$. This answers a question of A.~Valette.
\end{abstract}

\ead{piotr.nowak@vanderbilt.edu}


\end{frontmatter}

Let $X$ be a normed space. An \emph{affine, isometric action} of a
group $\G$ on $X$ is defined as $\Psi(g)v=\pi(g)v+\gamma(g)$ for
$v\in X$, $g\in \G$, where $\pi$ is a unitary (i.e. linear
isometric) representation of $\G$ on $X$ and $\gamma:\G\to X$
satisfies the cocycle identity with respect to $\pi$, i.e.
$\gamma(gh)=\pi(g)\gamma(h)+\gamma(g)$. The action is
\emph{proper} if $\lim_{g\to\infty} \Vert \Psi(g)v\Vert=\infty$
for every $v\in X$. This is equivalent to $\lim_{g\to\infty} \Vert
\gamma(g)\Vert=\infty$. One can express this idea in the language
of coarse geometry by saying that every orbit map is a coarse
embedding.

The following definition is due to Gromov.

\begin{defn}[{\cite[6.A.III]{gromov}}]
A second countable, locally compact group is said to be
\emph{a-T-menable} (has the \emph{Haagerup approximation
property}) if it admits a proper, affine, isometric action on a
separable Hilbert space $\Hb$.
\end{defn}

A-T-menability was designed as a strong opposite of Kazhdan's
property (T). We recall here a geometric characterization of
property (T) known as the Delorme-Guichardet Theorem, for a
detailed account of the subject see
\cite{bekka-delaharpe-valette}.

\begin{defn}
A second countable, locally compact group $\G$ has \emph{Kazhdan's
Property (T)} if and only if every affine isometric action of $\G$
on a Hilbert space has a fixed point.
\end{defn}
As suggested in the definition, a-T-menability turned out to be
equivalent to the Haagerup property (this was proved in
\cite{bekka-etal}), which arose in the study of approximation
properties of operator algebras and has application to harmonic
analysis. There are many other characterizations of a-T-menability,
in particular Gromov showed \cite[7.A]{gromov} that it is equivalent
to existence of a proper isometric action on the (either real or
complex) infinite dimensional hyperbolic space.

Recently N.~Brown and E.~Guentner \cite{brown-guentner} proved
that every discrete group admits a proper, affine and isometric
action on an $\ell_2$-direct sum $\left(\sum \ell_{p_n}\right)_2$,
for some sequence $\left\{p_n\right\}$ satisfying
${p_n}\longrightarrow\infty$. Since there are discrete groups
which are not a-T-menable, i.e. groups which are Kazhdan (T), an
existence of a proper, affine, isometric action on a reflexive
Banach space does not in general imply a-T-menability. Also
results of G.~Yu show that property (T) groups may admit proper,
affine, isometric actions on the spaces $\ell_p$ for $p>2$
\cite{yu-actions}. We also refer the reader to the recent article
\cite{bader-furman-gelander-monod} for a thorough study of similar
questions in the context of property (T).

What we are interested in is to find Banach spaces actions on
which imply or characterize a-T-menability. The motivation comes
from a question of A.Valette, who in \cite[Section
7.4.2]{cherix-etal} asked whether there are geometric
characterizations of a-T-menability other than through actions on
infinite-dimensional hyperbolic spaces. We prove the following

\begin{thm}\label{theorem1}
For a second countable, locally compact group $\G$ the following
conditions are equivalent:\begin{enumerate}\item $\G$ is
a-T-menable\\ \item $\G$ admits a proper, affine, isometric action
on the Banach space $L_p[0,1]$ for some $1<p<2$\\ \item $\G$
admits a proper, affine, isometric action on the Banach space
$L_p[0,1]$ for all $1<p<2$ \end{enumerate}
\end{thm}

Note that the results in \cite{brown-guentner,yu-actions} show
that Theorem \ref{theorem1} cannot be extended to $p>2$ or to the
class of reflexive or uniformly convex Banach spaces.

We also want to mention a problem raised in \cite[6.D$_3$]{gromov}
by Gromov: for a given group $\G$ find all such $p\ge 1$ for which
$\G$ admits a proper, affine, isometric action on $\ell_p$. Our
methods give some partial information on possible answers to this
question, namely Proposition \ref{prop.p.then.aT} states that only
a-T-menable groups may admit such actions on $\ell_p$ for $0<
p<2$.

A-T-menability is an important property in studying the
Baum-Connes Conjecture. N.~Higson and G.~Kasparov showed
\cite{higson-kasparov} that every discrete a-T-menable group
satisfies the Baum-Connes Conjecture with arbitrary coefficients.

\begin{ack}
The result of this paper arose from my work on the M.Sc. thesis
under Professor Henryk Toru\'nczyk at the University of Warsaw. I
would like to thank him for guidance. I am also grateful to
Norbert Riedel for many fruitful discussions, to Alain Valette for
helpful remarks and to the referee for improving the exposition of
this note.
\end{ack}

\section{Proofs}
 We will use the fact that a-T-menability can be
characterized in terms of existence of certain conditionally
negative definite functions, which we define now.

By a kernel on a set $X$ we mean a symmetric function $K:X\times
X\to\R$.

\begin{defn}

A kernel $K$ is said to be \emph{conditionally negative definite}
if $$\sum K(x_i,x_j)c_ic_j\le 0$$ for all $n\in\N$ and
$x_1,...,x_n\in X$, $c_1,...,c_n\in \R$ such that $\sum c_i=0$.

A function $\psi\colon \G\to \R$ on a metric group $\G$,
satisfying $\psi(g)=\psi(g^{-1})$ is said to be conditionally
negative definite if $K(g,h)=\psi(gh^{-1})$ is a conditionally
negative definite kernel.
\end{defn}

 It is easy to check that
if $(\Hb,\Vert\cdot\Vert)$ is a Hilbert space then the kernel
$K(x,y)=\Vert x-y\Vert^2$ is conditionally negative definite.

The following characterization is due to M.E.B.~Bekka,
P.-A.~Cherix and A.~Valette.
\begin{thm}[\cite{bekka-etal}]\label{a-T-menable.through.negative}
A second countable, locally compact group $\G$ is a-T-menable if
and only if there exists a continuous, conditionally negative
definite function $\psi:\G\to\R_+$ satisfying
$\lim_{g\to\infty}\psi(g)=\infty$.
\end{thm}

To prove Theorem \ref{theorem1} we also need the following lemmas
concerning conditionally negative definite functions and kernels
on $L_p$-spaces. These facts where proved by Schoenberg
\cite{schoenberg2}, for a further discussion see e.g.
\cite[Chapter 8]{bl}.

\begin{lem}\label{power.of.negative}
Let $K$ be a conditionally negative definite kernel on $X$ and
$K(x,y)\ge 0$ for all $x,y\in X$. Then the kernel $K^{\alpha}$ is
conditionally negative definite for any $0<\alpha<1$.
\end{lem}
\begin{pf*}{Proof.}
Let $K$ be a conditionally negative definite kernel.
 Then for every $t\ge 0$ the kernel $1-e^{-tK}\ge 0$ is also
conditionally negative definite and we have
$$\int\limits_0^{\infty}{\left(1-e^{-tK}\right)d\mu(t)}\ge 0$$
for every positive measure $\mu$ on $[0,\infty)$. For every $x>0$
and $0<\alpha<1$ the following formula holds
$$x^{\alpha}=c_{\alpha}\int\limits_0^{\infty}{\left(1-e^{-tx}\right)t^{-\alpha-1} dt},$$
where $c_{\alpha}$ is some positive constant. Thus $K^{\alpha}$ is
also a conditionally negative definite kernel for every
$0<\alpha<1$. \qed \end{pf*}

\begin{lem}\label{metric.on.Lp}
The function $\Vert x\Vert^p$ is conditionally negative definite
on $L_p(\mu)$ when $0< p\le 2$.
\end{lem}
\begin{pf*}{Proof.}
The kernel $\vert x-y\vert^2$ is conditionally negative definite
on the real line (as a square of the metric on a Hilbert space).
By Lemma \ref{power.of.negative}, for any $0<p\le 2$ the kernel
$\vert x-y\vert^p$ is also conditionally negative definite on
$\R$, i.e.,
$$\sum \vert x_i-x_j\vert^p c_i c_j\le 0$$
for every such $p$, all $x_1,...,x_n\in \R$ and $c_1,...,c_n\in\R$
such that $\sum c_i=0$. Integrate the above inequality with
respect to the measure $\mu$ to establish the proof. \qed
\end{pf*}

It follows from the lemmas that the norm on $L_p(\mu)$ is a
conditionally negative definite function, provided $1\le p\le2$.

To state the next proposition we define a more general notion of a
proper action, it is necessary when talking about the spaces
$L_p(\mu)$ for $p< 1$ which are not normable metric vector spaces.
Thus, if $X$ is just a metric space we call an isometric action of
$\G$ on $X$ proper if the set $\lbrace\,g\in\G\,\vert
g\mathcal{U}\cap\mathcal{U}\,\rbrace$ is finite for any bounded
set $\mathcal{U}\subset X$. For normed spaces this is consistent
with the definitions stated in the introduction.

\begin{prop}\label{prop.p.then.aT}
If a second countable, locally compact group $\G$ admits a proper,
affine, isometric action on a space $L_p(\mu)$ for some $0<p<2$
then $\G$ is a-T-menable.
\end{prop}
\begin{pf*}{Proof.}
Given a proper, affine, isometric $\G$-action on $L_p[0,1]$
consider the function $\psi:\G\to\R$, $\psi(g)=\Vert
\gamma(g)\Vert^p$, where $\gamma$ is the cocycle associated with
the action. Since the $p$-th power of the norm on $L_p[0,1]$ is a
conditionally negative definite function by Lemma
\ref{metric.on.Lp}, $\psi$ is a conditionally negative function on
$\G$. The considered $\G$-action is proper thus $\lim_{g\to\infty}
\Vert\gamma(g)\Vert^p=\infty$ and by Theorem
\ref{a-T-menable.through.negative}, $\G$ is a-T-menable.
\qed\end{pf*}
In particular only a-T-menable groups may admit
proper, affine isometric actions on the spaces $\ell_p$ for
$0<p<2$ (cf. \cite{yu-actions}).

\begin{pf*}{Proof of Theorem \ref{theorem1}.}
$(1) \Rightarrow (3)$. \ Let $G$ be a locally compact, second countable, a-T-menable group. Then by \cite[Theorem 2.2.2]{cherix-etal} there exists
a standard probability space $(X,\mu)$ and a measure preserving action of $G$ on $X$ such that
\begin{enumerate}
\item there exists a sequence of Borel sets $A_n\subseteq X$ such that $\mu(A_n)=\frac{1}{2}$ and
$\sup_{g\in B(e,n)} \mu(A_ng\triangle A_n)\le \frac{1}{2^n}$ ,
\item the action is strongly mixing, i.e. $\langle f,f\cdot g\rangle\to 0$ when $g\to\infty$ for
every $f\in L_2(X,\mu)$ such that $\int f d\mu=0$.
\end{enumerate}
Choose the sequence
Let $v_n(x)=1_{A_n}(x)-\frac{1}{2}\in L_2(X,\mu)$. Then  $\Vert v_n\Vert_2=\dfrac{1}{2}$ and
$$\int_{X}v_n(x)\ d\mu=0$$  
so by strong mixing, $$\Vert v_n-v_n\cdot g\Vert_2\to \sqrt{2}\Vert v_n\Vert_2,$$
when $g\to\infty$.
Also, for $g\in B(e,n)$ we have
$$\Vert v_n-v_n\cdot g\Vert_2=\mu(A_ng\,\triangle A_n)\le \dfrac{1}{2^n}$$
for all $g\in B(e,n)$.

Now given $p<2$ define
$$w_n(x)=\vert v_n(x)\vert^{2/p}\sign(v_n(x))\in L_p(X,\mu).$$
In other words, $w_n$ is a image of $v_n$ under the \emph{Mazur map}, which is a uniform
homeomorphism between unit balls of $L_p$-spaces, see \cite[Ch. 9.1]{bl}
for details and estimates. Moreover this map clearly commutes with the regular representation.
By the uniform continuity of the Mazur map and its inverse there exist constants  $C,\delta>0$ (which
depend only on $p$) such that the sequence $w_n$ satisfies
\begin{enumerate}
\item $\sup_{g\in B(e,n)}\Vert w_n\cdot s-w_n\Vert_p\le C/{2^n}$,
\item $\Vert w_n\cdot g-w_n\Vert_p \ge \delta$ for all $g\in G\setminus B(e,S_n)$ for some $S_n>0$, which depends on $n$ only (the sequence $\{S_n\}$ can be chosen to be increasing).
\end{enumerate}

This allows to construct a proper affine isometric action on $L_p(X,\mu)$ in a standard way.
Define  $b:G\to \left(\bigoplus_{n=1}^{\infty}L_p(X,\mu)\right)_p$
($p$ denotes the $L_p$-norm on the infinite direct sum) 
$$b(g)=\oplus_{n=1}^{\infty}\  \rho(g)w_n-w_n$$
where $\rho:G\to \mathrm{Iso}(L_p(X,\mu))$ is the right regular representation of $G$ on $X$.
Then $b$ is a cocycle for the representation $\oplus\, \rho$ by 
standard calculations (see e.g. \cite{bekka-etal}).

This way we obtain a proper isometric action on
$\left(\bigoplus_{n=1}^{\infty}L_p(X,\mu)\right)_p$ and the only thing left to notice is that
by construction in the proof of \cite[Theorem 2.2.2]{cherix-etal} the measure $\mu$ is non-atomic, thus by
the isometric classification of $L_p$-spaces, $L_p(X,\mu)$
is isometric to $L_p[0,1]$ and the $p$-sum of infinitely many of these spaces is again isometric to
$L_p[0,1]$. Thus $G$ admits a proper, affine, isometric action on $L_p[0,1]$.

$(3) \Rightarrow (2)$. This is obvious.

$(2) \Rightarrow (1)$. This implication is proved in Proposition
\ref{prop.p.then.aT} above. \qed\end{pf*}

Note that the above methods cannot be applied to other Banach
spaces. J.~Bretagnolle, D.~Dacuhna-Castelle and J.L.~Krivine
showed \cite{bdck} that the function $\Vert x\Vert^p$, $0<p\le2$,
is a conditionally negative definite kernel on a Banach space $X$
if and only if $X$ is isometric to a subspace of $L_p(\mu)$ for
some measure $\mu$. Together with Lemma \ref{power.of.negative}
this covers all powers $p\ge 1$.

\end{document}